\documentclass[a4paper,reqno,12pt]{amsart}
\usepackage{amsmath,amssymb,amsbsy,amsfonts,amsthm,latexsym,
            amsopn,amstext,amsxtra,euscript,amscd}
\usepackage{hyperref}   
\usepackage[margin=3cm]{geometry}

\begin{document}
\bibliographystyle{plain}
\newfont{\teneufm}{eufm10}
\newfont{\seveneufm}{eufm7}
\newfont{\fiveeufm}{eufm5}
%
%
\newfam\eufmfam
              \textfont\eufmfam=\teneufm \scriptfont\eufmfam=\seveneufm
              \scriptscriptfont\eufmfam=\fiveeufm
\def\bbbr{{\rm I\!R}}
\def\bbbm{{\rm I\!M}}
\def\bbbn{{\rm I\!N}}
\def\bbbf{{\rm I\!F}}
\def\bbbh{{\rm I\!H}}
\def\bbbk{{\rm I\!K}}
\def\bbbp{{\rm I\!P}}
\def\bbbone{{\mathchoice {\rm 1\mskip-4mu l} {\rm 1\mskip-4mu l}
{\rm 1\mskip-4.5mu l} {\rm 1\mskip-5mu l}}}
\def\bbbc{{\mathchoice {\setbox0=\hbox{$\displaystyle\rm C$}\hbox{\hbox
to0pt{\kern0.4\wd0\vrule height0.9\ht0\hss}\box0}}
{\setbox0=\hbox{$\textstyle\rm C$}\hbox{\hbox
to0pt{\kern0.4\wd0\vrule height0.9\ht0\hss}\box0}}
{\setbox0=\hbox{$\scriptstyle\rm C$}\hbox{\hbox
to0pt{\kern0.4\wd0\vrule height0.9\ht0\hss}\box0}}
{\setbox0=\hbox{$\scriptscriptstyle\rm C$}\hbox{\hbox
to0pt{\kern0.4\wd0\vrule height0.9\ht0\hss}\box0}}}}
\def\bbbq{{\mathchoice {\setbox0=\hbox{$\displaystyle\rm
Q$}\hbox{\raise
0.15\ht0\hbox to0pt{\kern0.4\wd0\vrule height0.8\ht0\hss}\box0}}
{\setbox0=\hbox{$\textstyle\rm Q$}\hbox{\raise
0.15\ht0\hbox to0pt{\kern0.4\wd0\vrule height0.8\ht0\hss}\box0}}
{\setbox0=\hbox{$\scriptstyle\rm Q$}\hbox{\raise
0.15\ht0\hbox to0pt{\kern0.4\wd0\vrule height0.7\ht0\hss}\box0}}
{\setbox0=\hbox{$\scriptscriptstyle\rm Q$}\hbox{\raise
0.15\ht0\hbox to0pt{\kern0.4\wd0\vrule height0.7\ht0\hss}\box0}}}}
\def\bbbt{{\mathchoice {\setbox0=\hbox{$\displaystyle\rm.
T$}\hbox{\hbox to0pt{\kern0.3\wd0\vrule height0.9\ht0\hss}\box0}}
{\setbox0=\hbox{$\textstyle\rm T$}\hbox{\hbox
to0pt{\kern0.3\wd0\vrule height0.9\ht0\hss}\box0}}
{\setbox0=\hbox{$\scriptstyle\rm T$}\hbox{\hbox
to0pt{\kern0.3\wd0\vrule height0.9\ht0\hss}\box0}}
{\setbox0=\hbox{$\scriptscriptstyle\rm T$}\hbox{\hbox
to0pt{\kern0.3\wd0\vrule height0.9\ht0\hss}\box0}}}}
\def\bbbs{{\mathchoice
{\setbox0=\hbox{$\displaystyle     \rm S$}\hbox{\raise0.5\ht0\hbox
to0pt{\kern0.35\wd0\vrule height0.45\ht0\hss}\hbox
to0pt{\kern0.55\wd0\vrule height0.5\ht0\hss}\box0}}
{\setbox0=\hbox{$\textstyle        \rm S$}\hbox{\raise0.5\ht0\hbox
to0pt{\kern0.35\wd0\vrule height0.45\ht0\hss}\hbox
to0pt{\kern0.55\wd0\vrule height0.5\ht0\hss}\box0}}
{\setbox0=\hbox{$\scriptstyle      \rm S$}\hbox{\raise0.5\ht0\hbox
to0pt{\kern0.35\wd0\vrule height0.45\ht0\hss}\raise0.05\ht0\hbox
to0pt{\kern0.5\wd0\vrule height0.45\ht0\hss}\box0}}
{\setbox0=\hbox{$\scriptscriptstyle\rm S$}\hbox{\raise0.5\ht0\hbox
to0pt{\kern0.4\wd0\vrule height0.45\ht0\hss}\raise0.05\ht0\hbox
to0pt{\kern0.55\wd0\vrule height0.45\ht0\hss}\box0}}}}
\def\bbbz{{\mathchoice {\hbox{$\sf\textstyle Z\kern-0.4em Z$}}
{\hbox{$\sf\textstyle Z\kern-0.4em Z$}}
{\hbox{$\sf\scriptstyle Z\kern-0.3em Z$}}
{\hbox{$\sf\scriptscriptstyle Z\kern-0.2em Z$}}}}
\def\ts{\thinspace}

\newtheorem{theorem}{Theorem}
\newtheorem{lemma}[theorem]{Lemma}
\newtheorem{claim}[theorem]{Claim}
\newtheorem{cor}[theorem]{Corollary}
\newtheorem{prop}[theorem]{Proposition}
\newtheorem{definition}[theorem]{Definition}
\newtheorem{remark}[theorem]{Remark}
\newtheorem{question}[theorem]{Open Question}
\newtheorem{example}[theorem]{Example}
\newtheorem{problem}[theorem]{Problem}

\def\qed{\ifmmode
\squareforqed\else{\unskip\nobreak\hfil
\penalty50\hskip1em\null\nobreak\hfil\squareforqed
\parfillskip=0pt\finalhyphendemerits=0\endgraf}\fi}

\def\squareforqed{\hbox{\rlap{$\sqcap$}$\sqcup$}}

\def \C {{\mathbb C}}
\def \F {{\mathbb F}}
\def \L {{\mathbb L}}
\def \K {{\mathbb K}}
\def \Q {{\mathbb Q}}
\def \Z {{\mathbb Z}}
\def\cA{{\mathcal A}}
\def\cB{{\mathcal B}}
\def\cC{{\mathcal C}}
\def\cD{{\mathcal D}}
\def\cE{{\mathcal E}}
\def\cF{{\mathcal F}}
\def\cG{{\mathcal G}}
\def\cH{{\mathcal H}}
\def\cI{{\mathcal I}}
\def\cJ{{\mathcal J}}
\def\cK{{\mathcal K}}
\def\cL{{\mathcal L}}
\def\cM{{\mathcal M}}
\def\cN{{\mathcal N}}
\def\cO{{\mathcal O}}
\def\cP{{\mathcal P}}
\def\cQ{{\mathcal Q}}
\def\cR{{\mathcal R}}
\def\cS{{\mathcal S}}
\def\cT{{\mathcal T}}
\def\cU{{\mathcal U}}
\def\cV{{\mathcal V}}
\def\cW{{\mathcal W}}
\def\cX{{\mathcal X}}
\def\cY{{\mathcal Y}}
\def\cZ{{\mathcal Z}}
\newcommand{\rmod}[1]{\: \mbox{mod}\: #1}

\def\tcN{\cN^\mathbf{c}}
\def\F{\mathbb F}
\def\Tr{\operatorname{Tr}}
\def\mand{\qquad \mbox{and} \qquad}
\renewcommand{\vec}[1]{\mathbf{#1}}
\def\eqref#1{(\ref{#1})}
\newcommand{\ignore}[1]{}
\hyphenation{re-pub-lished}
\parskip 1.5 mm
\def\lln{{\mathrm Lnln}}
\def\Res{\mathrm{Res}\,}
\def\F{{\bbbf}}
\def\Fp{\F_p}
\def\fp{\Fp^*}
\def\Fq{\F_q}
\def\ff{\F_2}
\def\ffn{\F_{2^n}}
\def\K{{\bbbk}}
\def \Z{{\bbbz}}
\def \N{{\bbbn}}
\def\Q{{\bbbq}}
\def \R{{\bbbr}}
\def \P{{\bbbp}}
\def\Zm{\Z_m}
\def \Um{{\mathcal U}_m}
\def \Bf{\frak B}
\def\Km{\cK_\mu}
\def\va {{\mathbf a}}
\def \vb {{\mathbf b}}
\def \vc {{\mathbf c}}
\def\vx{{\mathbf x}}
\def \vr {{\mathbf r}}
\def \vv {{\mathbf v}}
\def\vu{{\mathbf u}}
\def \vw{{\mathbf w}}
\def \vz {{\mathbfz}}
\def\\{\cr}
\def\({\left(}
\def\){\right)}
\def\fl#1{\left\lfloor#1\right\rfloor}
\def\rf#1{\left\lceil#1\right\rceil}
\def\flq#1{{\left\lfloor#1\right\rfloor}_q}
\def\flp#1{{\left\lfloor#1\right\rfloor}_p}
\def\flm#1{{\left\lfloor#1\right\rfloor}_m}
\def\Al{{\sl Alice}}
\def\Bob{{\sl Bob}}
\def\Or{{\mathcal O}}
\def\inv#1{\mbox{\rm{inv}}\,#1}
\def\invM#1{\mbox{\rm{inv}}_M\,#1}
\def\invp#1{\mbox{\rm{inv}}_p\,#1}
\def\Ln#1{\mbox{\rm{Ln}}\,#1}
\def \nd {\,|\hspace{-1.2mm}/\,}
\def\ord{\mu}
\def\E{\mathbf{E}}
\def\Cl{{\mathrm {Cl}}}
\def\epp{\mbox{\bf{e}}_{p-1}}
\def\ep{\mbox{\bf{e}}_p}
\def\eq{\mbox{\bf{e}}_q}
\def\bm{\bf{m}}
\newcommand{\floor}[1]{\lfloor {#1} \rfloor}
\newcommand{\comm}[1]{\marginpar{
\vskip-\baselineskip
\raggedright\footnotesize
\itshape\hrule\smallskip#1\par\smallskip\hrule}}
\def\rem{{\mathrm{\,rem\,}}}
\def\dist {{\mathrm{\,dist\,}}}
\def\etal{{\it et al.}}
\def\ie{{\it i.e. }}
\def\veps{{\varepsilon}}
\def\eps{{\eta}}
\def\ind#1{{\mathrm {ind}}\,#1}
               \def \MSB{{\mathrm{MSB}}}
\newcommand{\abs}[1]{\left| #1 \right|}

\title{ Power of Two as Sums of Three Pell Numbers}
%
\author{
{Jhon J. Bravo,}~~{\sc Bernadette~Faye}~~\and ~~ {\sc Florian~Luca}
}
\address{Departamento de Matemáticas, Universidad del Cauca, Calle 5 No. 4-70, Popay\'an, Colombia\newline
}
\email{jbravo@unicauca.edu.co}
\address{
African Institute for Mathematical Sciences(AIMS)\newline
Km 2 route de Joal\newline
BP 1418, Mbour, Senegal and\newline
School of Mathematics, University of the Witwatersrand \newline
Private Bag X3, Wits 2050, Johannesburg, South Africa
}
\email{bernadette@aims-senegal.org}

\address{
School of Mathematics, University of the Witwatersrand \newline
Private Bag X3, Wits 2050, Johannesburg, South Africa \newline
}
\email{Florian.Luca@wits.ac.za}

\begin{abstract} In this paper, we find all the  solutions of the Diophantine equation $P_\ell + P_m +P_n=2^a$,  in nonnegative integer variables $(n,m,\ell, a)$ where $P_k$ is the $k$-th term of the Pell sequence $\{P_n\}_{n\ge 0}$ given by $P_0=0$, $P_1=1$ and $P_{n+1}=2P_{n}+ P_{n-1}$ for all $n\geq 1$.
\end{abstract}

\maketitle



\section{Introduction}
\noindent The Pell sequence $\{P_n\}_{n\ge 0}$ is the binary reccurent sequence given by $P_0=0$, $P_1=1$ and $P_{n+1}=2P_{n}+ P_{n-1}$ for all $n\geq 0$. There are many papers in the literature dealing with Diophantine equations obtained by asking that
members of some fixed binary recurrence sequence be squares, factorials, triangular,  or belonging to some other interesting sequence of positive integers. 

For example, in $2008$, A. Peth\H o \cite{AP}  found all the perfect powers (of exponent larger than $1$) in the Pell sequence. His result is the following. 

\begin{theorem}[A. Peth\H o, \cite{AP}]
\label{th:1}
The only positive integer solutions  $(n,q,x)$ with $q\ge 2$ of the Diophantine equation 
$$P_n=x^q
$$ 
are $(n,q,x)=(1,q,1)$ and $(7,2,13)$. That is, the only perfect powers  of exponent larger than $1$ in the Pell numbers are 
$$
P_1=1 \quad \hbox{and}\quad P_7=13^2.
$$
\end{theorem}

The case $q=2$ had been treated earlier by Ljunggren \cite{LJ}. Peth\H o's result was rediscovered by J. H. E. Cohn \cite{Cohn}. 

In this paper, we study the following Diophantine equation: Find all nonnegative solutions $(\ell,m,n,a)$ of the equation
\begin{equation}
\label{eq:main}
P_\ell+P_m+P_n=2^a.
\end{equation}
There is already a vast literature on equations similar to \eqref{eq:main}. For example, putting for positive integers $a\ge 2$ and $n$, $s_a(n)$ for the sum of the base $a$ digits of $n$,   
Senge and Straus \cite{SS} showed that for each fixed $K$ and multiplicatively independent positive integers $a$ and $b$, the set $\{n: s_a(n)<K~{\text{\rm and}}~s_b(n)<K\}$ is finite. 
This was made effective by Stewart \cite{CL} using Baker's theory of lowers bounds for linear forms in logarithms of algebraic numbers (see also \cite{LuQu}). 
More concretely, the analogous equation \eqref{eq:main} when Pell numbers are replaced by Fibonacci numbers was solved in \cite{EJ3} (the special case when only two Fibonacci numbers are involved on the left has been solved 
earlier in \cite{EJ1}). Variants of this problem with  $k$-generalized Fibonacci numbers and Lucas numbers instead of Fibonacci numbers were studied in \cite{BGL16} and \cite{EJ2}, respectively. In \cite{BHL}, all Fibonacci numbers which are sums of three factorials were found, while in \cite{LuSi}, all  factorials which are sums of three Fibonacci numbers were found.  Repdigits which are sums of three Fibonacci numbers were found in \cite{LuRep}, while Fibonacci numbers which are sums of at most two repdigits were found in \cite{DL}. 

Our main result concerning \eqref{eq:1} is the following. 

\medskip

\begin{theorem}
\label{th:2}
The only solutions $(n,m,\ell,a)$ of the  Diophantine equation
\begin{equation}
\label{eq:1}
P_n+P_m+P_\ell = 2^a
\end{equation}
in integers $n\geq m \geq \ell \geq 0$ are in
$$
 \{(2,1,1,2),(3,2,1,3),(5,2,1,5),(6,5,5,7),(1,1,0,1),(2,2,0,2),(2,0,0,1),(1,0,0,0)\}.
$$
\end{theorem}

We use the method from \cite{LuRep}. 

\section{Preliminary results}

\noindent Let $(\alpha,\beta)=(1+{\sqrt{2}},1-{\sqrt{2}})$ be the roots of the characteristic equation $x^2-2x-1=0$ of the Pell sequence $\{P_n\}_{n\ge 0}$. The Binet formula for $P_n$ is
\begin{equation}
\label{eq:BinetP}
P_n= \frac{\alpha^n - \beta^n}{\alpha-\beta} \quad {\text{\rm for~ all}}\quad  n\ge 0.
\end{equation}
This implies easily that the inequalities
\begin{equation}
\label{eq:sizePn}
\alpha^{n-2}\le P_n\le  \alpha^{n-1}
\end{equation}
hold for all positive integers $n$.

Let $\{Q_n\}_{n\geq 0}$ be the companion Lucas sequence of the Pell sequence given by $Q_0=2$, $Q_1=2$ and $Q_{n+2}=2Q_{n+1}+ Q_{n}$ for all $n\ge 0$. For a prime $p$ and a nonzero integer $\delta$ let $\nu_p(\delta)$ be the exponent with which $p$ appears in the prime factorization of $\delta$.  The following result is well-known and easy to prove.
\begin{lemma}
\label{lem:orderof2}
The relations
\begin{itemize}
\item[(i)] $\nu_2(Q_n)=1$,
\item[(ii)] $\nu_2(P_n)=\nu_2(n)$
\end{itemize}
hold for all positive integers $n$.
\end{lemma}

The following result is an immediate consequence of Carmichael's primitive divisor theorem for Lucas sequences with real roots (see \cite{Car}).

\begin{lemma}
\label{lem:prim}
If $n\ge 13$, then $P_n$ has a prime factor $\ge n-1$.
\end{lemma}

We also need a Baker type lower bound for a nonzero linear form in logarithm of algebraic numbers. We choose to use the result of Matveev in \cite{MV}. Before proceeding further, we recall some basics notions from algebraic number theory.

Let $\eta$ be an algebraic number of degree $d$ over $\mathbb{Q}$ with minimal primitive polynomial over the integers
$$
f(X) = a_0 \prod_{i=1}^{d}(X-\eta^{(i)}) \in \mathbb{Z}[X],
$$
where the leading coefficient $a_0$ is positive and the $\eta^{(i)}$ are conjugates of $\eta$. The logarithmic
height of  $\eta$ is given by
$$
h(\eta) = \dfrac{1}{d}\left(\log a_0 + \sum_{i=1}^{d}\log\max\{|\eta^{(i)}|,1\}\right).
$$
The following properties of the logarithms height, which will be used in the next section without special reference, are also known:
\begin{itemize}
\item $h(\eta\pm \gamma)\leq h(\eta) + h(\gamma) + \log 2.$
\item $h(\eta\gamma^{\pm})\leq h(\eta) + h(\gamma).$
\item $h(\eta^{s})=|s|h(\eta).$
\end{itemize}

With these above notations, Matveev proved the following theorem (see also \cite{BMS}).
\begin{theorem}[Matveev \cite{MV}, Theorem 9.4 \cite{BMS}] 
\label{thm:Matveev}
Let $\K$ be a number field of degree $D$ over ${\mathbb Q}$, $\eta_1, \ldots, \eta_t$ be positive
real numbers of ${\mathbb K}$, and $b_1, \ldots,  b_t$ rational integers. Put
$$
\Lambda = \eta_1^{b_1} \cdots \eta_t^{b_t}-1
\qquad
\text{and}
\qquad
B \geq \max\{|b_1|, \ldots ,|b_t|\}.
$$
Let $A_i \geq \max\{Dh(\eta_i), |\log \eta_i|, 0.16\}$ be real numbers, for
$i = 1, \ldots, t.$
Then, assuming that $\Lambda \not = 0$, we have
$$
|\Lambda| > \exp(-1.4 \times 30^{t+3} \times t^{4.5} \times D^2(1 + \log D)(1 + \log B)A_1 \cdots A_t).
$$
\end{theorem}

In $1998$, Dujella and Peth\H o in \cite[Lemma 5$(a)$]{DP} gave a version of the reduction method originally proved by Baker and Davenport \cite{Baker-Davenport}.  We next present the following lemma from \cite{BL1} (see also \cite{BGL16}), which is an immediate variation of the result due to Dujella and Peth\H o from \cite{DP}, and is the key tool used to reduce the upper bound on the variable $n$. For a real number $x$ we put $\|x\|=\min\{|x-n|: n\in\mathbb{Z}\}$ for the distance from $x$ to the nearest integer.

\begin{lemma}
\label{reduce}
Let $M$ be a positive integer, let $p/q$ be a convergent of the continued fraction of the irrational $\gamma$ such that $q>6M$, and let $A,B,\mu$ be some real numbers with $A>0$ and $B>1$. Let $\epsilon:=||\mu q||-M||\gamma q||$. If $\epsilon >0$, then there is no solution to the inequality
$$
0<|u\gamma-v+\mu|<AB^{-w},
$$
in positive integers $u,v$ and $w$ with
$$
u\leq M \quad\text{and}\quad w\geq \frac{\log(Aq/\epsilon)}{\log B}.
$$
\end{lemma}

\section{Proof of Theorem \ref{th:2}}
\subsection{The case $\ell=0$}

If $\ell=m=0$, we then get that $P_n=2^a$. This implies that $n\le 12$ by Lemma \ref{lem:prim}.
If $\ell=0$ but $m>0$, we then get
\begin{equation}
\label{eq:tv}
P_n+P_m=2^a.
\end{equation}
Since $P_m$ and $P_n$ are positive, we get that $a>0$, so $P_n$ and $P_m$ have the same parity.  The left--hand side above factors as 
\begin{equation}
\label{eq:RelP1}
P_n+P_m= P_{(n+\delta m)/2}Q_{(n-\delta m)/2},
\end{equation}
where $\delta\in \{\pm 1\}$ is $1$ if $n\equiv m\pmod{4}$ and $-1$ otherwise, a fact easily checked. Thus, equation \eqref{eq:tv} becomes
$$
P_{(n+\delta m)/2}Q_{(n-\delta m)/2}=2^a.
$$
Lemmas \ref{lem:orderof2} and \ref{lem:prim} show that $(n-\delta m)/2\in \{0,1\}$ and $(n+\delta m)/2\le 12$, and all solutions can now be easily found. All in all, the case $\ell=0$ gives the last four solutions listed in the statement of Theorem \ref{th:2}.

\subsection{Bounding $n-m$  and $n-\ell$ in terms of $n$}

From now, we assume $n\geq m\geq \ell\geq 1$. First of all, if $n=m=\ell$, equation \eqref{eq:1} become $3P_n=2^a$ which is impossible. Thus, we assume from now that either $n> m$ or $m> \ell$. We next perform a computation showing to show that there are no others solutions to equation \eqref{eq:1} than those listed in Theorem \ref{th:2} in the range $1\leq \ell\leq m\leq n\leq 150.$ So, from now on we work under the assumption that $n>150$.

We find a relation between $a$ and $n$. Using equation \eqref{eq:1} and the right-hand side of inequality \eqref{eq:sizePn}, we get that
$$
2^a<\alpha^{n-1}+\alpha^{m-1} +\alpha^{\ell-1}< 2^{2n-2}( 1 + 2^{2(m-n)}+2^{2(\ell-n)})< 2^{2n+1}.
$$
where in the middle nequality we used the fact that $\alpha<2^2.$ Hence, we have that $a\leq 2n$.
We rewrite equation \eqref{eq:1} using \eqref{eq:BinetP} as
$$
\frac{\alpha^n}{2\sqrt{2}}-2^a=\frac{\beta^n}{2\sqrt{2}}-(P_m+ P_\ell).
$$
We take absolute values in both sides of the above relation  with the right-hand side of \eqref{eq:sizePn} obtaining
$$
\Big| \frac{\alpha^n}{2\sqrt{2}}-2^a\Big|\leq \frac{|\beta|^n}{2\sqrt{2}}+P_m+P_\ell< \frac{1}{2}+\(\alpha^m+\alpha^\ell\).
$$
Dividing both sides by $\alpha^n/(2\sqrt{2})$, we get 
\begin{equation}
\label{eq:2}
\Big| 1-2^{a+1}\cdot\alpha^{-n}\cdot\sqrt{2}\Big|< \frac{8}{\alpha^{n-m}}.
\end{equation}
We are in a situation to apply Matveev's result Theorem \ref{thm:Matveev} to the left--hand side of \eqref{eq:2}. The expression on the left-hand side of \eqref{eq:2} is nonzero, since this expression being zero means that $2^{a+1}=(\alpha^n/\sqrt{2})$, so $\alpha^{2n}\in {\mathbb Z}$ for some positive integer $n$, which is false. Hence, we take ${\mathbb K}:={\mathbb Q}({\sqrt{2}})$ for which $D=2$. We take 
$$
t:=3,\quad \eta_1:=2,\quad \eta_2:=\alpha,\quad \eta_3:=\sqrt{2},\quad b_1:=a+1,\quad b_2:=-n,\quad b_3:=1.
$$ 
So, we can take $A_1:=1.4$, $A_2:=0.9$ and $A_3:=0.7$. Finally we recall that $a\leq 2n$ and deduce that $\max\{|b_1|,|b_2|,|b_3|\}\leq 2n+1$, so we take $B:=2n+1$.
Theorem \ref{thm:Matveev} implies that a lower bound on the left-hand side of \eqref{eq:2} is
\begin{equation}
\label{eq:3}
\exp\left(-1.4\times 30^6\times 3^{4.5}\times  2^2\times (1+\log 2)(2\log n)\times 1.4\times 0.9\times 0.7\right).
\end{equation}
In the above inequality, we used $1+\log (2n+1)< 2\log n$, which holds in our range of $n$. Taking logarithms in inequality \eqref{eq:2} and comparing the resulting inequality with \eqref{eq:3}, we get that
\begin{equation}
\label{eq:4}
(n-m)\log \alpha < 1.8\times 10^{12}\log n.
\end{equation}
We now consider a second linear form in logarithms by rewriting equation \eqref{eq:1} in a different way. Using the Binet formula \eqref{eq:BinetP}, we get that
$$
\frac{\alpha^n}{2\sqrt{2}} +\frac{\alpha^m}{2\sqrt{2}}-2^a=\frac{\beta^n}{2\sqrt{2}} +\frac{\beta^m}{2\sqrt{2}}-P_\ell, 
$$
which implies
$$
\Big| \frac{\alpha^n}{2\sqrt{2}}(1+\alpha^{m-n})-2^a\Big|\leq \frac{|\beta|^n + |\beta|^m}{2\sqrt{2}}+P_\ell<\frac{1}{2}+\alpha^\ell.
$$
Dividing both sides of the above inequality by the first term of the left-hand side, we obtain
\begin{equation}\label{eq:5}
\Big| 1-2^{a+1}\cdot\alpha^{-n}\cdot\sqrt{2}(1+\alpha^{m-n})^{-1}\Big|< \frac{5}{\alpha^{n-\ell}}.
\end{equation}
We apply again Matveev Theorem \ref{thm:Matveev} with the same ${\mathbb K}$ as before. 
We take 
$$
t:=3,\quad  \eta_1:=2,\quad \eta_2:=\alpha,\quad \eta_3:=\sqrt{2}(1+\alpha^{m-n})^{-1},\quad  b_1:=a+1,\quad b_2:=-n,\quad b_3:=1.
$$ 
So, we can take $A_1:=1.4$, $A_2:=0.9$ and  $B:=2n+1$. We observe that the left-hand side of \eqref{eq:5} is not zero because otherwise we would get 
\begin{equation}
\label{eq:6}
2^{a+1}\sqrt{2}=\alpha^n(1+\alpha^{m-n})=\alpha^n+\alpha^m.
\end{equation}
By conjugating the above in $\mathbb{K}$ we get that 
\begin{equation}
\label{eq:7}
-2^{a+1}\sqrt{2}=\beta^n+\beta^m.
\end{equation}
Equations \eqref{eq:6} and \eqref{eq:7}, lead to
$$
\alpha^n<\alpha^n+\alpha^m=|\beta^n + \beta^m|\leq |\beta|^n+|\beta|^m<1,
$$
which is impossible. Now, let us have a look on the logarithmic height of $\eta_3$. Since,
$$
\eta_3=\sqrt{2}(1+\alpha^{m-n})^{-1}<\sqrt{2}\quad\hbox{and}\quad\eta_3^{-1}=\frac{1+\alpha^{m-n}}{\sqrt{2}}<\frac{2}{\sqrt{2}},
$$
we get that $|\log \eta_3|<1$. Furthermore, we notice that 
$$
h(\eta_3)\leq \log \sqrt{2}+|m-n|\(\frac{\log \alpha}{2}\)+\log 2=\log(2\sqrt{2})+(n-m)\(\frac{\log \alpha}{2}\).
$$
Thus, we can take $A_3:= 3+(n-m)\log\alpha>\max\{2h(\eta_3),|\log\eta_3|,0.16\}.$
As before, Theorem \ref{thm:Matveev} and \eqref{eq:5} imply that
\begin{equation}
\label{eq:8}
\exp\left(-2.45\times 10^{12}\times \log n\times (3+(n-m)\log\alpha)\right)<\frac{5}{\alpha^{n-\ell}}
\end{equation}
giving
\begin{equation}
\label{eq:9}
(n-\ell)\log\alpha < 2.5\times 10^{12}\times \log n\times (3+(n-m)\log\alpha).
\end{equation}
Inserting inequality \eqref{eq:4} into \eqref{eq:9}, we obtain
\begin{equation}
\label{eq:9bis}
(n-\ell)\log \alpha<5\times 10^{24}\log^2 n.
\end{equation}

\subsection{Bounding $n$}

We now use a third linear form in logarithms by rewriting equation \eqref{eq:1} in a different way. Using the Binet formula \eqref{eq:BinetP}, we get that

$$\frac{\alpha^n}{2\sqrt{2}} +\frac{\alpha^m}{2\sqrt{2}}+\frac{\alpha^\ell}{2\sqrt{2}}-2^a=\frac{\beta^n}{2\sqrt{2}} +\frac{\beta^m}{2\sqrt{2}}+\frac{\beta^\ell}{2\sqrt{2}}, $$
which implies

$$\Big| \frac{\alpha^n}{2\sqrt{2}}(1+\alpha^{m-n}+\alpha^{\ell-n})-2^a\Big|\leq \frac{|\beta|^n + |\beta|^m + |\beta|^\ell}{2\sqrt{2}}<\frac{1}{2}$$
for all $n> 150$ and $m\geq\ell\geq 1$. Dividing both sides of the above inequality by the first term of the left-hand side, we obtain

\begin{equation}
\label{eq:5bis}
\Big| 1-2^{a+1}\cdot\alpha^{-n}\cdot\sqrt{2}(1+\alpha^{m-n} + \alpha^{\ell-n})^{-1}\Big|< \frac{2}{\alpha^{n}}.
\end{equation}
As before, we use Matveev Theorem \ref{thm:Matveev} with the same  ${\mathbb K}$ as before and with 
$$
t:=3, ~ \eta_1:=2,~\eta_2:=\alpha,~ \eta_3:=\sqrt{2}(1+\alpha^{m-n}+\alpha^{\ell-n})^{-1},\quad b_1:=a+1,~ b_2:=-n, b_3:=1.
$$ 
As before we take $A_1:=1.4$, $A_2:=0.9$ and  $B:=2n+1$. It remains us to prove that the left-hand side of \eqref{eq:5bis} is not zero. Assuming the contrary, we would get 
\begin{equation}
\label{eq:6bis}
2^{a+1}\sqrt{2}=\alpha^n(1+\alpha^{m-n}+ \alpha^{\ell-n})=\alpha^n+\alpha^m + \alpha^{\ell}.
\end{equation}
Conjugating the above relation in $\mathbb{K}$ we get that 
\begin{equation}
\label{eq:7bis}
-2^{a+1}\sqrt{2}=\beta^n+\beta^m + \beta^{\ell}.
\end{equation}
Equations \eqref{eq:6bis} and \eqref{eq:7bis}, lead to
$$
\alpha^n<\alpha^n+\alpha^m +\alpha^\ell =|\beta^n + \beta^m+ \beta^{\ell}|\leq |\beta|^n+|\beta|^m + |\beta|^\ell<1
$$
which is impossible since $\alpha>2$. It remains to estimate the logarithmic height of $\eta_3$. Since,
$$
\eta_3=\sqrt{2}(1+\alpha^{m-n}+ \alpha^{\ell-n})^{-1}<\sqrt{2}\quad\hbox{and}\quad\eta_3^{-1}=\frac{1+\alpha^{m-n}+\alpha^{\ell-n}}{\sqrt{2}}<\frac{3}{\sqrt{2}},
$$
it follows that $|\log \eta_3|<1.$, Furthermore, we notice that 

\begin{eqnarray*}
h(\eta_3) &\leq &\log \sqrt{2}+|m-n|\(\frac{\log \alpha}{2}\)+ |\ell-n|\(\frac{\log \alpha}{2}\)+2\log 2 \\
&=&\log(4\sqrt{2})+(n-m)\(\frac{\log \alpha}{2}\) + (n-\ell)\(\frac{\log \alpha}{2}\).
\end{eqnarray*}
Thus, we can take 
$$
A_3:= 4+(n-m)\log\alpha+ (n-\ell)\log\alpha>\max\{2h(\eta_3),|\log\eta_3|,0.16\}.
$$
As before, Theorem \ref{thm:Matveev} and \eqref{eq:5bis} implies that
\begin{equation}\label{eq:8bis}
\exp\left(-2.45\times 10^{12}\times \log n\times (4+(n-m)\log\alpha+ (n-\ell)\log\alpha)\right)< \frac{2}{\alpha^{n}}
\end{equation}
which leads to
\begin{equation}
\label{eq:91bis}
n\log\alpha < 2.5\times 10^{12}\times \log n\times (4+(n-m)\log\alpha+ (n-\ell)\log\alpha).
\end{equation}
Inserting inequalities \eqref{eq:4} and \eqref{eq:9bis} into \eqref{eq:91bis} and performing the required computations, we obtain
\begin{equation}
\label{eq:92bis}
n<1.7\times 10^{37}\log^3 n,
\end{equation}
giving $n<1.7\times 10^{43}.$ We summarize the conclusion of this section as follows.
\medskip

\begin{lemma}
\label{lem:1}
If $(n,m,\ell, a)$ is a solution in positive integers of equation \eqref{eq:1}, with $n\ge m\ge \ell$, 
then 
$$
a< 2n+1<4\times 10^{43}.
$$
\end{lemma}

\subsection{Reducing the bound on $n$}

We use several times Lemma \ref{reduce} to reduce the bound for $n$. We return to \eqref{eq:2}. Put
$$
\Lambda_1:=(a+1)\log 2-n\log \alpha+\log\sqrt{2}.
$$
Then \eqref{eq:2} implies that 
\begin{equation}
\label{eq:10}
|1-e^{\Lambda_1}|<\frac{8}{\alpha^{n-m}}.
\end{equation}
Note that $\Lambda_1>0$ since 
$$
\frac{\alpha^n}{2\sqrt{2}}<P_n+1\leq P_n+P_m+P_\ell=2^{a}.
$$
Hence, using the fact that $1+x<e^x$ holds for all positive real numbers $x$, we get that
$$
0<\Lambda_1\leq e^{\Lambda_1}-1<\frac{8}{\alpha^{n-m}}.
$$
Dividing across by $\log\alpha$, we get
\begin{equation}
\label{eq:12}
0< (a+1)\left(\frac{\log 2}{\log\alpha}\right)-n+\left(\frac{\log\sqrt{2}}{\log\alpha}\right)<\frac{10}{\alpha^{n-m}}.
\end{equation}
We are now ready to apply Lemma \ref{reduce} with the obvious parameters
$$
\gamma:=\frac{\log 2}{\log \alpha},\quad \mu:=\frac{\log\sqrt{2}}{\log\alpha},\quad A:=10,\quad B:=\alpha.
$$
It is easy to see that $\gamma$ is irrationnal. We can take $M:=4\times 10^{43}$. Applying Lemma \ref{reduce} and performing the calculations with $q_{91}>6M$ and $\epsilon:=||\mu q_{91}||-M||\gamma q_{91}||>0$, we get that if $(n,m,\ell,a)$ is a solution to equation \eqref{eq:1}, then $n-m\in[0,130]$. 
We now work with inequality \eqref{eq:5} to obtain an upper bound on $n-\ell$. We put 
$$
\Lambda_2:=(a+1)\log 2-n\log \alpha+\log g(n-m),
$$
where we put $g(x):=\sqrt{2}(1+\alpha^{-x})^{-1}$. Then \eqref{eq:5} implies that 
\begin{equation}
\label{eq:11bis}
|1-e^{\Lambda_2}|<\frac{5}{\alpha^{n-\ell}}.
\end{equation}
Using the Binet formula of the Pell sequence with \eqref{eq:1}, one can show that $\Lambda_2> 0$ since 
$$
\frac{\alpha^n}{2\sqrt{2}}+\frac{\alpha^m}{2\sqrt{2}}<P_n+P_m+1\leq P_n+P_m+P_\ell=2^{a}.
$$
From this and \eqref{eq:11bis} we get
$$
0<\Lambda_2<\frac{5}{\alpha^{n-\ell}}.
$$
Replacing $\Lambda_2$ in the above inequality by its formula and arguing as in \eqref{eq:12}, we get that 
\begin{equation}
\label{eq:13}
0< (a+1)\left(\frac{\log 2}{\log\alpha}\right)-n+\frac{\log g(n-m)}{\log\alpha}<\frac{6}{\alpha^{n-\ell}}.
\end{equation}
Here, we take $M:=4\times 10^{43}$ and as we explain before, we apply Lemma \ref{reduce} to inequality \eqref{eq:13} for all possible choices of $n-m \in [0,130]$, except when $n-m=1,2$. Computing all the possible cases with suitable values for the parameter $q$, we find that if $(n,m,\ell,a)$ is a solution of \eqref{eq:1}, with $n-m\neq 1,2$, then $n-\ell\leq 140$.
For the special cases where $n-m=1,2$,  we have that 
\begin{equation*}
\label{eq:RelP}
\frac{\log g(x)}{\log\alpha}=\left\{ \begin{matrix}
0 & {\text{if}} & x=1;\vspace{0.3cm}\\
1-\frac{\log 2}{\log\alpha} & {\text{if}} & x=2.
\end{matrix}\right.
\end{equation*}
Thus, we cannot apply Lemma \ref{reduce}, because the value for the parameter $\varepsilon$ is always $\leq 0.$ Thus, in these cases, the reduction algorithm is not useful. However, we can see that if $n-m=1,2$, then the resulting inequality from \eqref{eq:13} has the shape $0 <|x\gamma-y|<6/\alpha^{n-\ell}$ with $\gamma$ being an irrational number and $x,y \in\mathbb{Z}$. So, we can appeal to the known properties of the convergents of the continued fractions to obtain a nontrivial lower bound for $|x\gamma-y|$. This gives us an upper bound for $n-\ell$.  Let's see the details. 
When $n-m=1$, $\log g(n-m)/\log\alpha=0$ and we get from \eqref{eq:13} that 
\begin{equation}
\label{eq:15}
0< (a+1)\gamma-n<\frac{6}{\alpha^{n-\ell}} \quad\hbox{where}\quad \gamma:=\frac{\log 2}{\log\alpha}.
\end{equation}
Let $[a_0,a_1,a_2,\ldots]=[0, 1, 3, 1, 2,\ldots]$ be the continued fraction expression of the above $\gamma$ and let $p_k/q_k$ be the its $k$⁻th convergents. Recall that $a+1< 4\times 10^{43}.$ A quick computation with Mathematica shows that 
$$q_{87}<4\times 10^{43}<q_{88}.$$
Furthermore $a_M:=\max\{a_i: i=1\ldots,88\}=100.$ Then, from the properties of the continued fractions, inequality \eqref{eq:15} becomes
$$
\frac{1}{(a_M+2)(a+1)}< (a+1)\gamma-n<\frac{6}{\alpha^{n-\ell}}
$$
which yields to 
\begin{equation}
\label{exp1}
\alpha^{n-\ell}<6\cdot 102\cdot 4\times 10^{43}.
\end{equation}
Thus, $n-\ell<122.$ The same argument as before gives that $n-\ell<122$ in the case when $n-m=2.$ Therefore, $n-\ell\leq 140$ always holds.
Finally, in order to obtain a better upper bound on $n$, we use again inequality \eqref{eq:5bis} where we put 
$$
\Lambda_3:=(a+1)\log 2-n\log \alpha+\log \phi(n-m,n-\ell),
$$
with $\phi(x_1,x_2):=\sqrt{2}(1+\alpha^{-x_1}+\alpha^{-x_2})^{-1}$. Then \eqref{eq:5bis} implies that 
\begin{equation}
\label{eq:11}
|1-e^{\Lambda_3}|<\frac{2}{\alpha^{n}}.
\end{equation}
We observe that $\Lambda_3\neq 0$. We now analyze the cases $\Lambda_3>0$ and $\Lambda_3<0$. If $\Lambda_3>0,$ then 
$$
0<\Lambda_3<\frac{2}{\alpha^n}.
$$
Suppose now that $\Lambda_3<0.$ Since $2/\alpha^n<1/2$ for $n>150$, from \eqref{eq:11}, we get that $|e^{\Lambda_3}-1|<1/2$, therefore $e^{|\Lambda_3|}<2.$ Since $\Lambda_3<0$, we have that 
$$
0<|\Lambda_3|\leq e^{|\Lambda_3|}-1=e^{|\Lambda_3|}|e^{\Lambda_3}-1|<\frac{4}{\alpha^n}.
$$
Thus, we get in both cases that 
$$
0<|\Lambda_3|<\frac{4}{\alpha^n}.
$$
Replacing $\Lambda_3$ in the above inequality by its formula and arguing as in \eqref{eq:12}, we get that 
\begin{equation}
\label{eq:16}
0< \left|(a+1)\left(\frac{\log 2}{\log\alpha}\right)-n+\left(\frac{\log \phi(n-m,n-\ell)}{\log\alpha}\right)\right|<\frac{5}{\alpha^{n}}.
\end{equation}
Here, we take $M:=4\times 10^{43}$ and as we explained before, we apply Lemma \ref{reduce} to inequality \eqref{eq:16} for all possible choices of $n-m \in [0,130]$ and $n-\ell\in [0,140]$. With the help of Mathematica, we find that if $(n, m, \ell, a)$ is a possible solution of the equation \eqref{eq:1}, then $n < 150$, contradicting our assumption that $n>150$.  This finishes the proof of the theorem. 

\section*{Acknowledgments} 
\noindent J.~J.~B.~was supported in part by Project VRI ID 3744 (Universidad del Cauca). B. F. thanks AIMS for the AASRG. Her work in this project were carried out with financial support from the government of Canada's International Development Research Centre (IDRC) and whithin the framework of the AIMS Research for Africa Project.

\end{document}